\documentclass[letterpaper, 10pt, conference]{ieeeconf}

\IEEEoverridecommandlockouts
\usepackage{cite}
\usepackage[pdftex]{graphicx}
\graphicspath{{graphics/}} % declare the path(s) where your graphic files are
\DeclareGraphicsExtensions{.eps,.jpg,.pdf,.png}

\usepackage{balance}
\usepackage{tikz}
\usetikzlibrary{fadings}
\usetikzlibrary{patterns}
\usetikzlibrary{shadows.blur}
\usetikzlibrary{shapes}

\usepackage{amsfonts}
\usepackage{amsmath}
\usepackage{amssymb}
\usepackage{amsthm}
\allowdisplaybreaks
\usepackage{mathtools} % Additional symbols and fixes
\usepackage{upgreek}
\newcommand{\argmin}{\arg\!\min}

\DeclareMathOperator{\prox}{\operatorname{prox}}
\DeclareMathOperator{\sgn}{\operatorname{sgn}}

\newtheoremstyle{plain}
  {\topsep}   % ABOVESPACE
  {\topsep}   % BELOWSPACE
  {}          % BODYFONT
  {0pt}       % INDENT (empty value is the same as 0pt)
  {\bfseries} % HEADFONT
  {}          % HEADPUNCT
  {5pt plus 1pt minus 1pt} % HEADSPACE
  {}          % CUSTOM-HEAD-SPEC
\newtheoremstyle{plainit}
  {\topsep}   % ABOVESPACE
  {\topsep}   % BELOWSPACE
  {\itshape}  % BODYFONT
  {0pt}       % INDENT (empty value is the same as 0pt)
  {\bfseries} % HEADFONT
  {}          % HEADPUNCT
  {5pt plus 1pt minus 1pt} % HEADSPACE
  {}          % CUSTOM-HEAD-SPEC

\theoremstyle{plain}

\theoremstyle{plainit}

\newtheorem{theorem}{Theorem}
\newtheorem{lemma}{Lemma}

\newtheorem{remark}{Remark}

\usepackage{algorithmicx}
\usepackage{algorithm}
\usepackage{algpseudocode}

\usepackage{subcaption}

\usepackage{epstopdf}
\usepackage{url}
\newcommand{\loadhyperref}{
\usepackage[
  colorlinks=true,        % links are colored
  citecolor=black,        % color of cite links
  linkcolor=black,        % color of hyperref links
  menucolor=black,        % color of Acrobat Reader menu buttons
  urlcolor=black,         % color of urls \url{}
  bookmarks=true, hyperindex=true, breaklinks=true
]{hyperref}
}
\loadhyperref

\usepackage[capitalize]{cleveref}
\crefname{assumption}{Assumption}{Assumptions}
\crefname{corollary}{Corollary}{Corollaries}
\crefname{theorem}{Theorem}{Theorems}
\crefname{lemma}{Lemma}{Lemmas}

\newcommand{\yss}{y_{\textup{ss}}}

\newcommand{\rev}[1]{\textcolor{black}{#1}}

\definecolor{lightgreen}{RGB}{199, 237, 204}
\def\arxivVersion{}

\begin{document}

\title{\LARGE \bf Robust Feedback Optimization with Model Uncertainty: \\ A Regularization Approach \\
\thanks{
\textsuperscript{*}Equal contribution. \textsuperscript{1}Automatic Control Laboratory, ETH Z{\"u}rich, Switzerland. Email: \{wechan, zhiyhe, kmoffat, bsaverio, dorfler\}@ethz.ch. \textsuperscript{2}Max Planck Institute for Intelligent Systems, T{\"u}bingen, Germany. Email: \{zhiyu.he, michael.muehlebach\}@tuebingen.mpg.de. This work was supported by the Max Planck ETH Center for Learning Systems, the SNSF via NCCR Automation (grant agreement 51NF40 80545), and the German Research Foundation.}
}

\author{Winnie~Chan\textsuperscript{1,*}, Zhiyu He\textsuperscript{1,2,*}, Keith Moffat\textsuperscript{1}, Saverio Bolognani\textsuperscript{1}, Michael Muehlebach\textsuperscript{2}, and Florian D{\"o}rfler\textsuperscript{1}
}

\maketitle
\thispagestyle{empty}
\pagestyle{empty}

\begin{abstract}
Feedback optimization optimizes the steady state of a dynamical system by implementing optimization iterations in closed loop with the plant. It relies on online measurements and limited model information, namely, the input-output sensitivity. In practice, various issues, including inaccurate modeling, lack of observation, or changing conditions, can lead to sensitivity mismatches, causing closed-loop sub-optimality or even instability. To handle such uncertainties, we pursue robust feedback optimization, where we optimize the closed-loop performance against all possible sensitivities lying in specific uncertainty sets. We provide tractable reformulations for the corresponding min-max problems via regularizations and characterize the online closed-loop performance through the tracking error in case of time-varying optimal solutions. Simulations on a distribution grid illustrate the effectiveness of our robust feedback optimization controller in addressing sensitivity mismatches in a non-stationary environment.
% with similar interpretations as ridge and lasso regression
% finite-time tracking error as a function of the system property and the rate of change of optimal solutions. 

% Feedback optimization (FO) drives systems to optimal steady states without detailed model information or extensive numerical computations. However, inaccuracies in input-output sensitivities used for gradient-based updates can destabilize FO, leading to suboptimal performance or system failure. We propose robust FO, which incorporates sensitivity uncertainties to provide coupled optimality and stability guarantees for linear systems. The robust objective is reformulated into two min-max problems corresponding to ridge and Lasso regularization, for which we employ efficient controller designs based on projected (proximal) gradient descent. Our methodology adapts to time-varying environments and achieves dynamic regret proportional to the optimal solution's path length. Empirical results demonstrate superior performance to standard FO under model mismatch, particularly in high-dimensional scenarios. We validate our approach on the IEEE 123-bus system, showing effective control under varying operating regimes and topology change, even when standard FO fails.
\end{abstract}

% Note that keywords are not normally used for peer review papers.
% \begin{IEEEkeywords}
% feedback optimization, robust control, stability, model mismatch, power systems.
% \end{IEEEkeywords}

% !TEX root = ..\main.tex
\section{Introduction}
Modern engineering systems are increasingly complex, large-scale, and variable, as seen in power grids, supply chains, and recommender systems. Achieving optimal steady-state operation of these systems is both critical and challenging. In this regard, numerical optimization pipelines operate in an open-loop manner, whereby solutions are found based on an explicit formulation of the input-output map of the system and knowledge of disturbances. However, the reliance on accurate models poses restrictions and renders these pipelines unfavorable in complex environments.

Feedback optimization is an emerging paradigm for steady-state optimization of a dynamical system \cite{simonetto2020time,hauswirth2021optimization,krishnamoorthy2022real}. At the heart of feedback optimization is the interconnection between an optimization-based controller and a physical system. This closed-loop approach shares a similar spirit to extremum seeking\cite{ariyur2003real}, modifier adaptation\cite{marchetti2009modifier}, and real-time iterations\cite{diehl2005real}. Nonetheless, feedback optimization effectively handles high-dimensional objectives and coupling constraints, adapts to non-stationary conditions, and entails less computational effort (see review in \cite{hauswirth2021optimization}).

% At the interface of control, optimization, and learning, feedback optimization cross-fertilizes insights from different fields and offers flexibility, e.g., in handling constraints, fulfilling adaptation.

Thanks to the iterative structure that incorporates real-time measurements and performance objectives, feedback optimization enjoys closed-loop stability\cite{colombino_towards_2019}, optimality\cite{hauswirth_timescale_2021,lawrence_linear-convex_2021}, constraint satisfaction\cite{bianchin2021time}, and online adaptation\cite{colombino_online_2020,belgioioso2022online,ospina2022feedback,cothren2022online,ortmann2023deployment}. However, these salient properties rely on limited model information, i.e., the input-output sensitivity of a system. This requirement follows from using the chain rule to construct gradients in iterative updates. In practice, different issues can render the sensitivity inaccurate or elusive, e.g., corrupted data, lack of measurements, or changing conditions. As we will show in \cref{subsec:example_divergence}, such sensitivity errors can accumulate in the closed loop and cause significant sub-optimality or even divergence.

% Feedback optimization brings salient properties, e.g., closed-loop stability, optimality, and adaptation, when an accurate sensitivity matrix of the system is available \cite{simonetto2020time,hauswirth2021optimization}.

Many approaches have been developed to address inexact sensitivities in feedback optimization. A major stream leverages model-free iterations, where controllers entirely bypass sensitivities. Such model-free operations are typically enabled by derivative-free optimization, including Bayesian\cite{simonetto2021personalized,krishnamoorthy2023model,xu_violation-aware_2024} and zeroth-order optimization\cite{poveda2017robust,chen2020model,tang2023zeroth,he2022model,chen2025continuous}. However, controllers based on Bayesian optimization tend to be computationally expensive for high-dimensional problems, whereas zeroth-order feedback optimization brings increased sample complexity. Therefore, it is desirable to incorporate structural, albeit inexact, sensitivity information into controller iterations rather than discard it altogether.

There are two primary solutions to handle model uncertainty without resorting to model-free iterations: adaptation and robustness. In the context of feedback optimization, adaptive schemes leverage offline or online data to refine knowledge of sensitivities, thereby facilitating closed-loop convergence. Examples include learning sensitivity via least squares\cite{picallo2021adaptive,dominguez2023online} or stochastic approximation\cite{agarwal2023model}, as well as constructing behavioral representations of sensitivity from input-output data\cite{bianchin2021online}. However, adaptive strategies impose additional requirements for data, computation, and estimation. Restrictions arise in scenarios involving high-dimensional systems and limited computational power, where sensitivity estimation can be challenging.

In this paper, we consider \emph{robust} feedback optimization, where the closed-loop performance is optimized given the worst-case realization of the sensitivity in some uncertainty sets. This is formalized as a min-max problem for which tractable reformulations via regularization are further provided. Our robust feedback optimization controllers feature provable convergence guarantees for time-varying problems with changing disturbances and references. Compared to the above adaptive schemes, our controllers only leverage an inexact sensitivity and hence are easy to implement. In contrast to related robust strategies in learning\cite{xu2010robust,bertsimas2011theory} and data-driven control\cite{huang2023robust}, we tackle a more demanding setting wherein model uncertainty is intertwined with both system dynamics and controller iterations. Our main contributions are as follows.

\begin{itemize}
    \item We formulate robust feedback optimization by addressing structured uncertainties in sensitivities. We provide tractable reformulations via regularization and build connections with lasso and ridge regression.
    % To our best knowledge, this is the first work in FO dealing explicitly with robustness to major model mismatch. We propose computationally-tractable robust objective formulations to address these uncertainties, and expound on workable controller design.

    \item We present online robust feedback optimization controllers that address two types of sensitivity uncertainty sets. We establish closed-loop convergence by characterizing errors in tracking trajectories of time-varying optimal solutions.
    % We present both theoretical and numerical analyses that characterize the performance of our robust control strategies. Premising both are mathematical and simulated examples of instability under the standard, non-robust FO regime.

    \item Through a numerical experiment of voltage regulation in a distribution grid, we demonstrate that the proposed controllers preserve voltage stability while prescribing less curtailment and reactive power control, even with inaccurate sensitivities.
    % We deploy robust FO in power system simulations of a distribution grid and demonstrate its practical utility in preserving voltage stability while prescribing less curtailment and reactive power control than normal FO controllers.
\end{itemize}

The rest of this paper is organized as follows. \cref{sec:problemFormulation} motivates and presents the problem setup. \cref{sec:design} provides tractable reformulations and our robust feedback optimization controllers. The closed-loop performance guarantee is established in \cref{sec:guarantee}, followed by numerical evaluations on a distribution grid in \cref{sec:simulation}. Finally, \cref{sec:conclusion} concludes the article and discusses future directions.

% \subsection*{Organization - to rewrite}
% Dummy text here.
% This paper is structured as such: \autoref{chap:problemFormulation} sets the stage for FO in linear dynamical systems and shows how it can be destabilized, before presenting robust formulations of the initial problem under two norm-defined uncertainty sets. Subsequently, \autoref{chap:controllers} details how controllers can carry out optimization with these robust objectives, while \autoref{chap:perfAnalysis} follows with a three-part theoretical analysis of robust FO's performance in terms of its coupled errors from input suboptimality and distance from steady state, regret order, and regularization effects. In \autoref{chap:experiments}, a battery of synthetic simulations are run to empirically compare standard and robust FO, as well as examine the sensitivity of robust FO's regularization parameter. We then bring our proposed control scheme to a power systems setting of voltage regulation amidst a high penetration of distributed generation in \autoref{chap:powerSystems} and illustrate its use in three operational scenarios. The last section ensues with a summary of our report and a discussion of robust FO's outlook.

% !TEX root = ..\main.tex
\section{Background and Problem Formulation}\label{sec:problemFormulation}
\subsection{Preliminaries}
We consider the following dynamical system
\begin{equation} \label{eqn:plant-dynamics}
    \begin{split}
        x_{k+1} &= Ax_k + Bu_k + d_{x, k}, \\
        y_k &= Cx_k + d_{y, k},
    \end{split}
\end{equation}
where $x_k \in \mathbb{R}^n$, $u_k \in \mathbb{R}^m$, $y_k \in \mathbb{R}^p$, $d_{x,k} \in \mathbb{R}^n$, and $d_{y,k} \in \mathbb{R}^p$ denote the state, input, output, exogenous disturbance, and measurement noise at time $k$, respectively. Further, $A \in \mathbb{R}^{n \times n}$, $B \in \mathbb{R}^{n \times m}$, and $C \in \mathbb{R}^{p \times n}$ are system matrices. 
% Throughout the paper,
We focus on a stable system, i.e., the spectral radius $\rho(A)$ of $A$ in \eqref{eqn:plant-dynamics} lies in $(0,1)$. In practice, this condition also holds if this system is prestabilized \rev{by state feedback controllers}. Given fixed inputs and disturbances (i.e., $u_k = u, d_{x,k} = d_x, d_{y,k} = d_y, \forall k \in \mathbb{N}$), system~\eqref{eqn:plant-dynamics} admits a unique steady-state output
\begin{equation}\label{eq:ss_output_linear_sys}
    \begin{split}
        \yss(u,d) &= Hu + d, \\
        H &\triangleq C(I\!-\!A)^{-1}B, \quad d \triangleq C(I\!-\!A)^{-1}d_x + d_y.
    \end{split}
\end{equation}
In \eqref{eq:ss_output_linear_sys}, $H \in \mathbb{R}^{p \times m}$ is the sensitivity matrix of system \eqref{eqn:plant-dynamics}.

A performance objective characterizing the input-output performance of system~\eqref{eqn:plant-dynamics} at each time $k \in \mathbb{N}$ is
\begin{equation}\label{eq:obj_quadratic}
	\begin{split}
    	\Phi_k(u; d_k, r_k) &= \|u\|^2_R + \lambda \|\yss(u, d_k) - r_k\|^2_Q \\
    			  &= \|u\|^2_R + \lambda \|H u + d_k - r_k\|^2_Q,
    \end{split}
\end{equation}
where $R \in \mathbb{R}^{m \times m}$ and $Q \in \mathbb{R}^{p \times p}$ are positive semidefinite matrices, $\|u\|_R \!=\! \sqrt{u^\top R u}$ and $\|y\|_Q \!=\! \sqrt{y^\top Q y}$ denote weighted norms, $\lambda \geq 0$ is a weight parameter, and $r_k \in \mathbb{R}^p$ is the reference at time $k$. Further, $\yss(u, d_k) = Hu + d_k$ is the steady-state output associated with the input $u$ and the disturbance $d_k \triangleq C(I\!-\!A)^{-1}d_{x,k} + d_{y,k}$ at time $k$. The function \eqref{eq:obj_quadratic} penalizes the input cost and the difference between the steady-state output and the reference.

To optimize \eqref{eq:obj_quadratic}, numerical solvers require an explicit knowledge of the map $\yss$ as per \eqref{eq:ss_output_linear_sys} with an accurate value of the disturbance $d$, which can be restrictive in applications. In contrast, feedback optimization leverages real-time output measurements and the limited model information, namely, the sensitivity matrix $H$, thus steering system~\eqref{eqn:plant-dynamics} to optimal operating conditions \cite{simonetto2020time,hauswirth2021optimization}.

\subsection{Example: Detrimental Effects of Inexact Sensitivities}\label{subsec:example_divergence}
Many practical issues including data inadequacy and varying conditions cause model uncertainty, i.e., sensitivity errors\cite{hauswirth2021optimization}. We present a motivating example to show how such errors invalidate feedback optimization by inducing closed-loop sub-optimality or instability. While this example is synthetic, we observe similar phenomena in realistic power grid simulations (see \cref{sec:simulation}).  

\begin{figure}[t!]
    \centering
    \begin{subfigure}[t]{0.49\columnwidth}
        \centering
        \includegraphics[width=\linewidth]{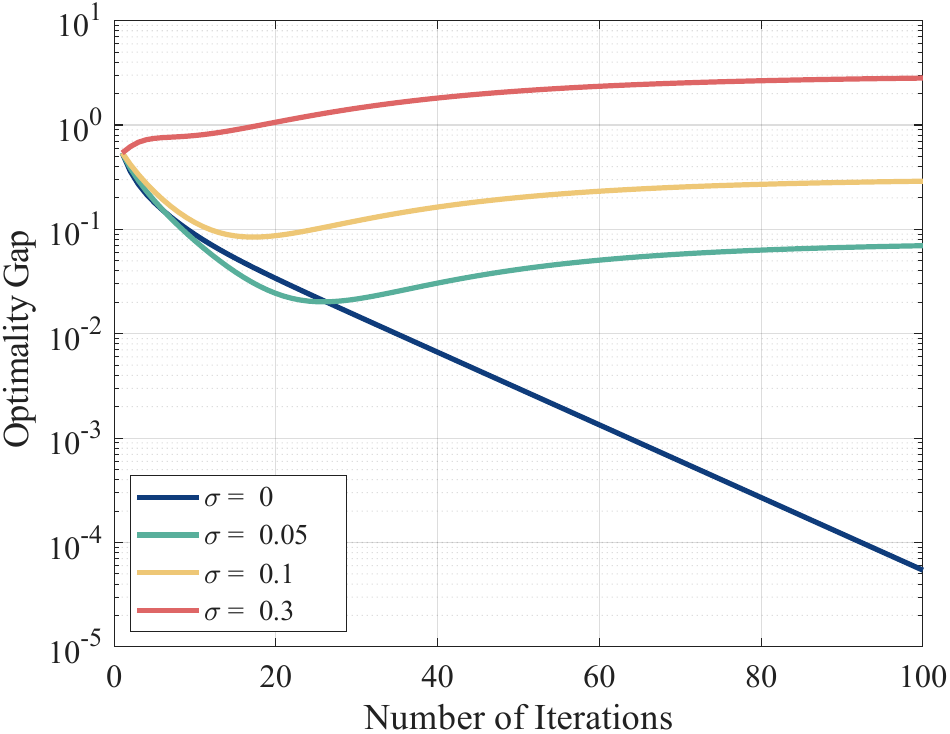}
        \caption{Use inexact sensitivities with varying perturbations.}
        \label{fig:inexact_sens_increase_noise}
    \end{subfigure}%
    \hfill
    \begin{subfigure}[t]{0.49\columnwidth}
        \centering
        \includegraphics[width=\linewidth]{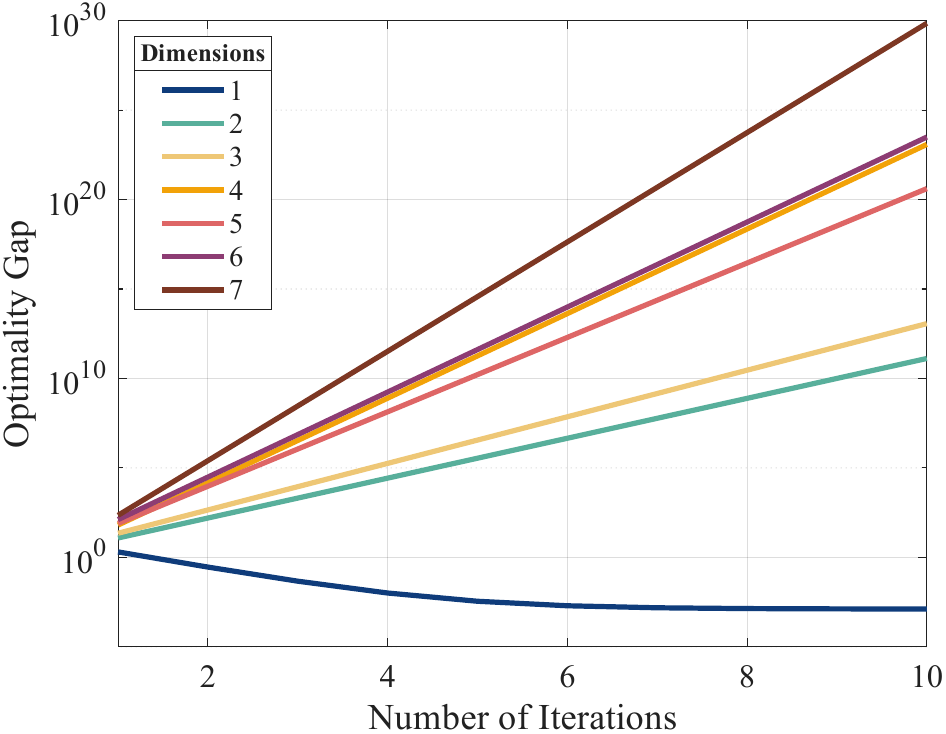}
        \caption{Use inexact sensitivities of different sizes.}
        \label{fig:inexact_sens_increase_dim}
    \end{subfigure}
    \caption{Closed-loop performance when the controller \eqref{eq:fo_inexact_sens} with inexact sensitivities is interconnected with the system \eqref{eq:ss_output_linear_sys}.}
    \label{fig:example_inexact_sens}
\end{figure}

We consider a system abstracted by the steady-state map \eqref{eq:ss_output_linear_sys} with fixed disturbances. We generate inexact sensitivities $\hat{H}$ in the following two fashions.
\begin{itemize}
    \item We fix the size of $H$ (i.e., $H \in \mathbb{R}^{3 \times 3}$) and add constant perturbations of different magnitudes. Specifically, $\hat{H} = H + \sigma \Delta_H$, where $\sigma \geq 0$, and the elements of $\Delta_H$ follow uniform distributions.
    % on $[-1, 1]$
    \item We vary the order of a square sensitivity $H$ from $1$ to $7$ and add perturbation noise with fixed norms, i.e., $\hat{H} = H + \Delta_H$. The square of each element of $\Delta_H$ satisfies a Dirichlet distribution, ensuring $\lVert \Delta_H \rVert_F = 1$.
\end{itemize}
To optimize \eqref{eq:obj_quadratic} \rev{with $r_k=0$}, consider the following feedback optimization controller using an inexact $\hat{H}$
% and acting on the system \eqref{eqn:plant-dynamics}:
\begin{equation}\label{eq:fo_inexact_sens}
    u_{k+1} = u_k - 2\eta \left(R u_k + \lambda \hat{H}^{\top} Q y_k \right),
\end{equation}
where $\eta > 0$ is the step size. The update \eqref{eq:fo_inexact_sens} follows a gradient descent iteration given the objective \eqref{eq:obj_quadratic}, replacing the steady-state output $H u_k + d$ by the real-time output measurement $y_k$ of \eqref{eqn:plant-dynamics}. \rev{We calculate the optimal value of \eqref{eq:obj_quadratic} offline through \textsf{fmincon} in \textsf{MATLAB} with the exact $H$ and $d$.} \cref{fig:example_inexact_sens} illustrates the closed-loop optimality gap when the controller \eqref{eq:fo_inexact_sens} is applied to the system \eqref{eqn:plant-dynamics}. We observe from \cref{fig:inexact_sens_increase_noise} that larger errors in sensitivities cause increased sub-optimality. Furthermore, \cref{fig:inexact_sens_increase_dim} demonstrates that when $\eta$ is fixed, the detrimental effect becomes more pronounced as the problem dimension grows.

\subsection{Problem Formulation}\label{subsec:formulation}
Motivated by the above observations, we pursue robust feedback optimization, where we optimize a worst-case performance objective given any realization of sensitivity lying in uncertainty sets. In practice, we can obtain through prior knowledge or identification\cite{ortmann2023deployment,hauswirth2021optimization,he2024gray} an inexact sensitivity $\hat{H}$, which differs from the true sensitivity $H$ of \eqref{eqn:plant-dynamics} by $\Delta_H$, i.e., $\hat{H} + \Delta_H = H$. In view of $\hat{H}$ and the uncertainty $\Delta_H$, our robust formulation is
\begin{equation} \label{eqn:robobj}
    \min_{u \in \mathbb{R}^m} \max_{\Delta_H \in \mathcal{D}} \quad \lVert u \rVert ^2_R + \lambda \lVert (\hat{H} + \Delta_H) u + d_k - r_k \rVert ^2_Q,
\end{equation}
% \zhmargin{We may need unconstrained problems}
where $\mathcal{D} \subset \mathbb{R}^{p \times m}$ is the uncertainty set wherein  $\Delta_H$ lies, $d_k = C(I\!-\!A)^{-1}d_{x,k} + d_{y,k}$ aggregates the disturbances $d_{x,k}$ and $d_{y,k}$, and $r_k$ is the reference at time $k$. Different from \eqref{eq:obj_quadratic}, in \eqref{eqn:robobj} we robustify the steady-state specification of system \eqref{eqn:plant-dynamics} against the sensitivity uncertainty $\Delta_H$. Essentially, \eqref{eqn:robobj} implies minimizing the steady-state input-output
performance for the worst-case sensitivity realization. 
% Similar to \eqref{eq:obj_quadratic}, the intuition of \eqref{eqn:robobj} is still to optimize the steady-state input-output performance of \eqref{eqn:plant-dynamics}. However, we address the uncertainty in $\hat{H}$ and aim to minimize the performance objective for the worst-case realization of sensitivity.
% $\mathcal{U} \subset \mathbb{R}^m$ is the input constraint set
% $\hat{\updelta}_k = \hat{d}_k - r_k$
% The circumflex accent $\hat{\cdot}$ thus indicates uncertainty in the parameters it describes. In the following, robust FO is examined under two prescribed settings for the uncertainty set.
We examine the following types of uncertainty sets.
\begin{itemize}
	\item \emph{Generalized uncertainties} described by
	\begin{equation}\label{eq:uncertainty_set_gen}
	    \mathcal{D}_{\text{gen}} \triangleq \left\{\Delta_H \big| \lVert \lambda^{\frac{1}{2}} Q^{\frac{1}{2}} \Delta_H \rVert _F \leq \varrho_{\text{gen}} \right\},
	\end{equation}
    where $\varrho_{\text{gen}} \geq 0$, and $\|\cdot\|_F$ denotes the Frobenius norm.

	\item \emph{Uncorrelated column-wise uncertainties} of the form
	\begin{align}
	    \mathcal{D}_{\text{col}} \triangleq \big\{\Delta_H \big| & \lVert (\lambda^{\frac{1}{2}} Q^{\frac{1}{2}} \Delta_H)_i \rVert \leq (\varrho_{\text{col}})_i, \notag \\
        &\forall i \in \{1, \cdots, m \} \big\}, \label{eq:uncertainty_set_col}
	\end{align}
    where $(\lambda^{\frac{1}{2}} Q^{\frac{1}{2}} \Delta_H)_i$ denotes the $i$-th column of the matrix $\lambda^{\frac{1}{2}} Q^{\frac{1}{2}} \Delta_H$ and $(\varrho_{\text{col}})_i$ denotes the $i$-th element of the vector $(\varrho_{\text{col}}) \in \mathbb{R}^m$, with $(\varrho_{\text{col}})_i \geq 0$.
    % where $(\lambda^{\frac{1}{2}} Q^{\frac{1}{2}} \Delta_H)_i$ and $(\varrho_{\text{col}})_i$ denote the $i$-th column and the $i$-th element of the matrix $\lambda^{\frac{1}{2}} Q^{\frac{1}{2}} \Delta_H$ and the vector $(\varrho_{\text{col}}) \in \mathbb{R}^m$, respectively, and $(\varrho_{\text{col}})_i \geq 0$.
\end{itemize}
In the above sets, $\mathcal{D}_{\text{gen}}$ poses a bounded-norm restriction on the uncertainty $\Delta_H$. In contrast, $\mathcal{D}_{\text{col}}$ bounds the norm of each column of $\Delta_H$, which is useful when different levels of confidence exist regarding how each component of $u$ affects the output $y$. Both types are common in the robust optimization literature \cite{bertsimas2011theory,xu2010robust,huang2023robust}.

% Although the entries in each row of $S^{\frac{1}{2}} \Delta_{H}$ are correlated to one \textbf{}another, entries in different rows and columns are not; thus, the \textit{column-wise} uncertainties are uncorrelated.
While problem \eqref{eqn:robobj} is unconstrained, we will discuss strategies to handle input and output constraints at the end of \cref{subsec:controller_design}. We consider quadratic objectives in \eqref{eqn:robobj} to highlight intuition and facilitate the presentation of robust strategies. Promising extensions to handle general objectives can be built on modern advances in robust optimization \cite{ben2008selected}.

\section{Robust Feedback Optimization} \label{sec:design}
\subsection{Tractable Reformulations} \label{subsec:reformulation}
We provide tractable reformulations of problem~\eqref{eqn:robobj}, thereby facilitating the design of robust feedback optimization controllers. Let $u_k^* \in \mathbb{R}^m$ be the optimal point of problem~\eqref{eqn:robobj} at time $k$. \rev{The objective function of \eqref{eqn:robobj} can be compactly written as
\begin{equation} \label{eqn:compact-robobj-formal}
    \Bigg\lVert \bigg( \underbrace{\begin{bmatrix} R^{\frac{1}{2}} \\ \lambda^{\frac{1}{2}} Q^{\frac{1}{2}}\hat{H} \end{bmatrix}}_{\triangleq M \in \mathbb{R}^{(m\!+\!p) \times m}} + \underbrace{\begin{bmatrix} 0 \\ \lambda^{\frac{1}{2}} Q^{\frac{1}{2}} \Delta_H \end{bmatrix}}_{\triangleq \Delta_M \in \mathbb{R}^{(m\!+\!p) \times m}} \bigg) u + \underbrace{\begin{bmatrix} 0 \\ d_k - r_k \end{bmatrix}}_{\triangleq \varepsilon_k \in \mathbb{R}^{m\!+\!p}} \Bigg\rVert ^2.
\end{equation}
We analyze two cases involving the uncertainty sets discussed in \cref{subsec:formulation}.}

\smallskip
\noindent\emph{$\bullet$ The case with generalized uncertainties}
\smallskip

For problem~\eqref{eqn:robobj} with the uncertainty set $\mathcal{D} = \mathcal{D}_{\text{gen}}$ (see \eqref{eq:uncertainty_set_gen}), the reformulated problem is
\begin{equation} \label{eqn:obj-l2}
    \min_{u \in \mathbb{R}^m} ~ \Phi_{k,\ell_2}(u) \triangleq \lVert u \rVert ^2_R + \lambda \lVert \hat{H} u + d_k - r_k \rVert ^2_Q + \rho_{\textup{gen}} \lVert u \rVert ^2.
\end{equation}
\rev{where the regularizer satisfies
\begin{equation}\label{eq:condition_l2}
    \rho_{\text{gen}} \!=\! 
    \begin{cases}
        \varrho_{\text{gen}} \lVert M u_k^* \!+\! \varepsilon_k \rVert / \|u_k^*\|,\! & \text{if } M u_k^* \!+\! \varepsilon_k \neq 0, \\
        \varrho_{\text{gen}} / \|u_k^*\|, & \text{otherwise},
    \end{cases}
\end{equation}
and $M$ and $\varepsilon_k$ are given in \eqref{eqn:compact-robobj-formal}.}

\smallskip
\noindent\emph{$\bullet$ The case with uncorrelated column-wise uncertainties}
\smallskip

The reformulated problem associated with \eqref{eqn:robobj} involving $\mathcal{D} = \mathcal{D}_{\text{col}}$ (see \eqref{eq:uncertainty_set_col}) is
\begin{equation} \label{eqn:obj-l1}
    \min_{u \in \mathbb{R}^m} ~ \Phi_{k,\ell_1}(u) := \lVert u \rVert ^2_R + \lambda \lVert \hat{H} u + d_k - r_k \rVert ^2_Q + \rho_{\textup{col}}^{\top} \lvert u \rvert,
\end{equation}
where $|u| \in \mathbb{R}^m$ denotes the component-wise absolute value of $u$, and the regularizer $\rho_{\textup{col}} \in \mathbb{R}^m$ satisfies
\begin{equation}\label{eq:condition_l1}
    \rho_{\textup{col}} = 
    \begin{cases}
        2 \lVert M u_k^* + \varepsilon_k \rVert \varrho_\textup{col}, & \text{if } M u_k^* + \varepsilon_k \neq 0, \\
        \varrho_{\textup{col}}, & \text{otherwise}.
    \end{cases}
\end{equation}

The following theorem establishes that the above reformulated problems share the same optimal points as problem~\eqref{eqn:robobj}.

\begin{theorem}
    Problem~\eqref{eqn:robobj} with the uncertainty set $\mathcal{D} = \mathcal{D}_{\textup{gen}}$ in \eqref{eq:uncertainty_set_gen} and problem~\eqref{eqn:obj-l2} share the same optimal point.
    % A vector $u_*$ is a minimizer of \eqref{eqn:robobj} if and only if $u_*$ also minimizes
    Moreover, problem~\eqref{eqn:robobj} with $\mathcal{D} = \mathcal{D}_{\textup{col}}$ in \eqref{eq:uncertainty_set_col} and problem~\eqref{eqn:obj-l1} attain the same optimal point.
    % i.e., $\frac{\lVert \Phi u_* + \varepsilon \rVert}{\sqrt{\lVert u_* \rVert ^2 + 1}}$
\end{theorem}

% For generalized uncertainties:
% \begin{lemma} \label{lem:robust-l2-fo}
%     A vector $u_*$ is a minimizer of
%     \begin{equation} \label{eqn:obj-l2-og}
%         \min_{u \in \mathcal{U}} \max_{(\Delta_H, \Delta_d) \in \mathcal{D}_{\text{gen}}} \quad \lVert u \rVert ^2_R + \lambda \lVert (\hat{H} + \Delta_H) u + \hat{\updelta} + \Delta_d \rVert ^2_S
%     \end{equation}
%     if and only if $u_*$ also minimizes
%     \begin{equation} \label{eqn:obj-l2}
%         \min_{u \in \mathcal{U}} \quad f_{l_2}(u) := \lVert u \rVert ^2_R + \lambda \lVert \hat{H} u + \hat{\updelta} \rVert ^2_S + \rho_{\textup{gen}} \lVert u \rVert ^2,
%     \end{equation}
%     where $\varrho_{\textup{gen}}$ and $\rho_{\textup{gen}} \in \mathbb{R}_{\geq 0}$ are equivalent up to a scaling constant of $\frac{\lVert \Phi u_* + \varepsilon \rVert}{\sqrt{\lVert u_* \rVert ^2 + 1}}$ when $\Phi u_* +     \varepsilon \neq 0$, and are identical otherwise.
% \end{lemma}

\begin{proof}
% \ifdefined\arxivVersion
We analyze the compact representation \eqref{eqn:compact-robobj-formal} of the objective of problem \eqref{eqn:robobj}.
% We first compactly express the sum of the two norms in (\ref{eqn:robobj}) as
% \begin{equation} \label{eqn:compact-robobj}
%     \lVert (M + \Delta_{M}) u + \rev{\bar{d}_k - \bar{r}_k} \rVert ^2,
% \end{equation}
\rev{In \eqref{eqn:compact-robobj-formal}, $M = \big[R^{\frac{1}{2}}; \lambda^{\frac{1}{2}} Q^{\frac{1}{2}} \hat{H}\big] \in \mathbb{R}^{(m+p) \times m}$, $\Delta_M = \big[0; \lambda^{\frac{1}{2}} Q^{\frac{1}{2}} \Delta_H\big] \in \mathbb{R}^{(m+p) \times m}$, and $\varepsilon_k = [0; d_k - r_k] \in \mathbb{R}^{m+p}$}. Since the norm is non-negative, optimizing \eqref{eqn:compact-robobj-formal} is equivalent to optimizing $\lVert (M + \Delta_{M}) u + \rev{\varepsilon_k}\rVert$. When $\mathcal{D} = \mathcal{D}_{\text{gen}}$, we perform an analysis similar to that in \cite[Theorem 3.1]{el_ghaoui_robust_1997} and obtain
% by making the following reformulations:
% \begin{align*}
%     \Phi &= [R^{\frac{1}{2}}; \lambda^{\frac{1}{2}} S^{\frac{1}{2}} \hat{H}] \in \mathbb{R}^{(m+p) \times m}\\
%     \Delta_{\Phi} &= [0; \lambda^{\frac{1}{2}} S^{\frac{1}{2}} \Delta_H] \in \mathbb{R}^{(m+p) \times m}\\
%     \varepsilon &= [0; \lambda^{\frac{1}{2}} S^{\frac{1}{2}} (\hat{d} - r)] \in \mathbb{R}^{(m+p) \times 1}\\
%     \Delta_{\varepsilon} &= [0; \lambda^{\frac{1}{2}} S^{\frac{1}{2}} \Delta_d] \in \mathbb{R}^{(m+p) \times 1}.
% \end{align*}
% Due to the non-negativity of the norm, it is monotonic with respect to its square, and any optimization over (\ref{eqn:compact-robobj}) is equivalent to (i.e., has identical solutions to) that over the norm itself. This property is used to apply the adapted result of \cite[Theorem 3.1]{el_ghaoui_robust_1997} for a $\varrho_{\text{gen}}$ of arbitrary value:
\begin{equation*}
    \max_{\lVert \Delta_{M} \rVert_F \leq \varrho_{\text{gen}}} \lVert (M \!+\! \Delta_{M}) u \!+\! \varepsilon_k \rVert
    =
    \|M u + \varepsilon_k\| + \varrho_{\text{gen}} \|u\|.
\end{equation*}
Moreover, the following two problems
\begin{align*}
    \min_{u \in \mathbb{R}^m} \lVert M u \!+\! \rev{\varepsilon_k} \rVert \!+\! \varrho_{\text{gen}} \|u\| ~ \text{and} ~
    \min_{u \in \mathbb{R}^m} \lVert M u \!+\! \rev{\varepsilon_k} \rVert ^2 \!+\! \rho_{\text{gen}} \lVert u \rVert ^2
    % \label{eqn:compact-robobj-reg}   
\end{align*}
% \begin{equation} \label{eqn:l2-equivalence}
%     \min_{u \in \mathcal{U}} \quad \lVert \Phi u + \varepsilon \rVert + \varrho_{\text{gen}} \sqrt{\lVert u \rVert ^2 + 1}
%     \quad \text{and} \quad 
%     \min_{u \in \mathcal{U}} \quad \lVert \Phi u + \varepsilon \rVert ^2 + \rho_{\text{gen}} \lVert u \rVert ^2
% \end{equation}
share the same optimal point $u_k^*$ if the condition \eqref{eq:condition_l2} holds. 
This result can be proved by comparing the optimality conditions of both problems and noting that $0$ is a subgradient of $\|u\|$ at $u=0$. \rev{We further expand the above objective and obtain the equivalent problem \eqref{eqn:obj-l2}}.

We proceed to analyze the case when $\mathcal{D} = \mathcal{D}_{\text{col}}$. Analogous to \cite[Theorem 1]{xu2010robust}, we obtain
\begin{equation*}
    \max_{\substack{\lVert (\Delta_{M})_i \rVert \leq (\varrho_{\textup{col}})_i \\ \forall i \in \{1, \dots, m\}}} \lVert (M \!+\! \Delta_{M}) u + \varepsilon_k \rVert
    \!=\!
    \lVert M u + \varepsilon_k \rVert + \varrho_{\textup{col}}^{\top} |u|.
\end{equation*}
Furthermore, the optimal points of the following problems 
\begin{equation*} %\label{eqn:l1-equivalence}
    \min_{u \in \mathbb{R}^m} ~ \lVert M u + \varepsilon_k \rVert + \varrho_{\textup{col}}^{\top} |u|
    \quad \text{and} \quad 
    \min_{u \in \mathbb{R}^m} ~ \lVert M u + \varepsilon_k \rVert ^2 + \rho_{\textup{col}}^{\top} |u|
\end{equation*}
coincide when the regularizer satisfies the condition \eqref{eq:condition_l1}. \rev{Further expansion of the above objective leads to \eqref{eqn:obj-l1}.}
% \begin{equation*}
%     \rho_H = 
%     \begin{cases}
%         2 \lVert M u^* + \varepsilon \rVert \varrho_H, & \text{if } M u^* + \varepsilon \neq 0 \\
%         \varrho_H, & \text{otherwise}. \qedhere
%     \end{cases}
% \end{equation*}
% \else
%     See our online report \cite{chan2025robust}.
% \fi
\end{proof}

The $\ell_2$-regularizer in \eqref{eqn:obj-l2} and the $\ell_1$-regularizer in \eqref{eqn:obj-l1} admit the same interpretation as those in classical ridge and lasso regression\cite{bertsimas2011theory}. In essence, these regularization terms penalize the magnitude of the input, helping to achieve closed-loop stability in the face of model uncertainty. The use of the $\ell_1$-regularizer also promotes sparsity of control inputs.

\begin{remark}
    While the expressions of $\rho_\textup{gen}$ and $\rho_{\textup{col}}$ involve $\|u_k^*\|$, this dependence arises from the quadratic objective in \eqref{eqn:robobj} and the related step of equivalent reformulation. In practice, for a variation of \eqref{eqn:robobj} with non-squared $\ell_2$-norms, the regularizers in the reformulated problems will only entail the uncertainty bounds $\varrho_{\text{gen}}$ and $\varrho_{\text{col}}$ but not $\|u_k^*\|$.
\end{remark}

% \zhmargin{discuss the reliance on $u_k^*$}

% The expression in (\ref{eqn:obj-l1-colmax}) can be interpreted as adding a $l_1$-regularizer $\lVert u \rVert _1$ to a least squares problem that minimizes the squared error of $\lVert \Phi u + \varepsilon \rVert ^2$, such that Lasso regression is performed (cf. \cite{bertsimas2011theory, buhlmann_statistics_2011}).

\subsection{Design of Robust Feedback Optimization Controllers} \label{subsec:controller_design}
Based on the reformulations in \cref{subsec:reformulation}, we present our online robust feedback optimization controllers. These controllers leverage an inexact sensitivity $\hat{H}$ and real-time output measurements of system \eqref{eqn:plant-dynamics}. They employ optimization-based iterations, thereby driving the system to operating points characterized by \eqref{eqn:obj-l2} or \eqref{eqn:obj-l1}.

For problem \eqref{eqn:obj-l2} corresponding to the case with generalized uncertainties, our robust feedback optimization controller employs the following gradient-based update
\begin{equation} \label{eqn:l2-step}
    u_{k+1} = u_k - 2\eta \left(R u_k \!+\! \lambda \hat{H}^{\top} Q (y_k \!-\! r_k) \!+\! \rho_{\text{gen}} u_k\right),
\end{equation}
where $\eta > 0$ is the step size. The update direction of the controller \eqref{eqn:l2-step} is related to the negative gradient of \eqref{eqn:obj-l2} at $u_k$. Further, \eqref{eqn:l2-step} uses the output measurement $y_k$ of the true system \eqref{eqn:plant-dynamics} as feedback.
% by using the real-time measurement-based derivative of $f_{l_2, k}(u)$.

Problem \eqref{eqn:obj-l1}, associated with the case with uncorrelated column-wise uncertainties, involves a nonsmooth regularizer $\|u\|$. Hence, building on proximal gradient descent, the proposed controller updates as follows:
\begin{equation} \label{eqn:l1-step}
    u_{k+1} = \prox_{\eta \rho_{\textup{col}}} \big(u_k \!-\! 2\eta (R u_k \!+\! \lambda \hat{H}^{\top} Q (y_k \!-\! r_k) \big),
\end{equation}
where $\eta > 0$ is the step size. In \eqref{eqn:l1-step}, $\prox_{\eta \rho_{\textup{col}}}(u) \triangleq \argmin_{u' \in \mathbb{R}^m} \eta \rho_{\textup{col}}^\top |u| + \frac{1}{2}\|u' - u\|^2$ denotes the proximal operator of $\eta \rho_{\textup{col}}^\top |u|$, i.e., element-wise soft thresholding $\sgn(u_i) \max\{|u_i| - \eta (\rho_{\textup{col}})_i, 0\}$, where $\sgn(\cdot)$ is the sign function. Similar to \eqref{eqn:l2-step}, this controller uses the real-time output $y_k$ of system~\eqref{eqn:plant-dynamics} and iteratively calculates new inputs.
% by using the real-time measurement-based derivative of $f_{l_1, k}(u) - \rho_H^{\top} \lvert u \rvert$.

We further discuss various extensions for the proposed controllers \eqref{eqn:l2-step} and \eqref{eqn:l1-step}. In practice, restrictions on the input due to actuation limits or economic conditions can often be represented as a constraint set $\mathcal{U} \subset \mathbb{R}^m$. In this regard, we can project $u_k$ generated by \eqref{eqn:l2-step} and \eqref{eqn:l1-step} back to $\mathcal{U}$, thereby satisfying constraint satisfaction at every time step. Should output constraints be imposed e.g. from safety requirements, we can augment the objectives in \eqref{eqn:obj-l2} and \eqref{eqn:obj-l1} with suitable penalty (e.g., quadratic or log-barrier) functions and incorporate the resulting derivative terms into the updates \eqref{eqn:l2-step} and \eqref{eqn:l1-step}, see also \cite[Section 2.4]{hauswirth2021optimization}. 

% To handle the input constraint set $\mathcal{U}$ in \eqref{eqn:robobj}, we further augment the above iterations with projection.

% \hzy{projection}

% !TEX root = ..\main.tex
\section{Performance Guarantee} \label{sec:guarantee}
We present the performance guarantee of the closed-loop interconnection between system \eqref{eqn:plant-dynamics} and our robust feedback optimization controller. A major challenge is that sensitivity uncertainty is interlaced with system dynamics and controller iterations, complicating convergence analysis. To address this challenge, we analyze the coupled evolution of the system \eqref{eqn:plant-dynamics} and the proposed controller, while characterizing the cumulative effects of sensitivity uncertainty. 

Recall that $u_k^*$ is the optimal point of problem~\eqref{eqn:robobj} at time $k$, and that $d_k \triangleq C(I\!-\!A)^{-1}d_{x,k} + d_{y,k}$ aggregates the disturbances. We consider the stable system \eqref{eqn:plant-dynamics}, i.e., $\rho(A) < 1$. Let $x_{\textup{ss},k} \in \mathbb{R}^n$ be the steady state of \eqref{eqn:plant-dynamics} induced by $u_{k}$ and $d_{x,k}$. In other words, $x_{\textup{ss},k} = Ax_{\textup{ss},k} + Bu_{k} + d_{x,k}$, implying $x_{\textup{ss},k} = (I-A)^{-1}(Bu_{k} + d_{x,k})$. For any given positive definite $\bar{Q} \in \mathbb{R}^{n \times n}$, there exists a unique positive definite $P \in \mathbb{R}^{n \times n}$ satisfying the Lyapunov equation $A^\top P A- P + \bar{Q} = 0$. Let $\|x\|_P \triangleq \sqrt{x^\top P x}$ be the weighted norm and $\lambda_{\max}(P)$ be the maximum eigenvalue of $P$. Our performance guarantee is as follows.

\begin{theorem} \label{thm:bounds}
    Let system \eqref{eqn:plant-dynamics} be stable. There exists $\eta^* > 0$ such that for any $\eta \in (0, \eta^*]$, the closed-loop interconnection between \eqref{eqn:plant-dynamics} and the controller \eqref{eqn:l2-step} or \eqref{eqn:l1-step} guarantees
    % Let Assumptions \ref{ass:uxy-inputs} to \ref{ass:Pmat} hold. Then the optimality and stability errors can be jointly bounded as
    \begin{align} \label{eqn:optstab-1}
        \left \lVert
        \begin{bmatrix} \|u_{k} \!-\! u_{k}^*\| \\ \lVert x_{k} \!-\! x_{\textup{ss},k} \rVert_P \end{bmatrix}
        \right \rVert
        \leq& r_1 (c_M)^k \left \lVert
        \begin{bmatrix} \|u_0 \!-\! u_0^*\| \\ \lVert x_0 \!-\! x_{\textup{ss},0} \rVert_P \end{bmatrix}
        \right \rVert \notag \\
        &+ r_2 \frac{c_M}{1 - c_M} \left \lVert \sup_{i \in [k]} q_i \right \rVert,
    \end{align}
    where $r_1, r_2 > 0$, and $c_M \in [0, 1)$. Moreover,
    \begin{equation*}
        q_k \triangleq \begin{bmatrix}
        \eta \bar{c}_1 \|H\!-\!\hat{H}\| \|u_{k}\| \!+\! \lVert u_{k+1}^* \!-\! u_k^* \rVert \\
        \eta c_3 \|H \!-\! \hat{H}\|\|u_k\| + c_4 \|d_{k+1} \!-\! d_{k}\| + \eta c_5
    \end{bmatrix},
    \end{equation*}
    where the constants are $\bar{c}_1 = 2\lambda \|\hat{H}^\top Q\|$, $c_3 = 2 \lambda \|(I\!-\!A)^{-1}B\| \lambda_{\max}(P) \|\hat{H}^\top Q\|$, $c_4 = \lambda_{\max}(P)\|(I\!-\!A)^{-1}\|$, and
    \begin{equation*}
        c_5 = \begin{cases} 0, & \textup{for } \eqref{eqn:l2-step}, \\ 2 \lambda_{\max}(P) \|(I-A)^{-1}B\| \|\rho_{\textup{col}}\|, & \textup{for } \eqref{eqn:l1-step}. \end{cases}
    \end{equation*}
\end{theorem}

\begin{proof}
\ifdefined\arxivVersion
The proof is provided in \cref{app:proof_thm}.
\else
See the full version of this paper online \cite{chan2025robust}. 
\fi
\end{proof}

In \cref{thm:bounds}, we characterize the closed-loop performance through the joint evolution of the distance $\|u_{k} - u_{k}^*\|$ to the optimal point $u_k^*$ and the distance $\|x_{k} - x_{\textup{ss},k}\|$ to the steady state $x_{\textup{ss},k}$. The upper bound \eqref{eqn:optstab-1} is in the flavor of input-to-state stability \cite{lasalle_stability_1986} and similar to \cite{belgioioso2022online,bianchin2021time,cothren2022online,ospina2022feedback}. In contrast to these works, we additionally characterize in \eqref{eqn:optstab-1} the cumulative effects of the given sensitivity uncertainty (i.e., $\|H - \hat{H}\|$) and the regularizer corresponding to the uncertainty set (i.e., $\rho_{\textup{col}}$). The effect of the initial conditions $u_0$ and $x_0$ vanishes exponentially fast, because $c_M \in (0,1)$. The asymptotic error is proportional to the shifts of optimal solutions $u_k^*$ and disturbances $d_k$, as well as the sensitivity uncertainty, i.e., $\|H - \hat{H}\|$. The influence arising from this uncertainty can be tuned via the step size, see the terms in $q_k$. It is possible to further establish upper bounds on the distance of the output to the optimal steady-state output through the Lipschitz property of the dynamics \eqref{eqn:plant-dynamics}.
\begin{figure*}[htb!]
    \centering
    \begin{subfigure}[t]{0.32\textwidth}
        \centering
        \includegraphics[width=\linewidth]{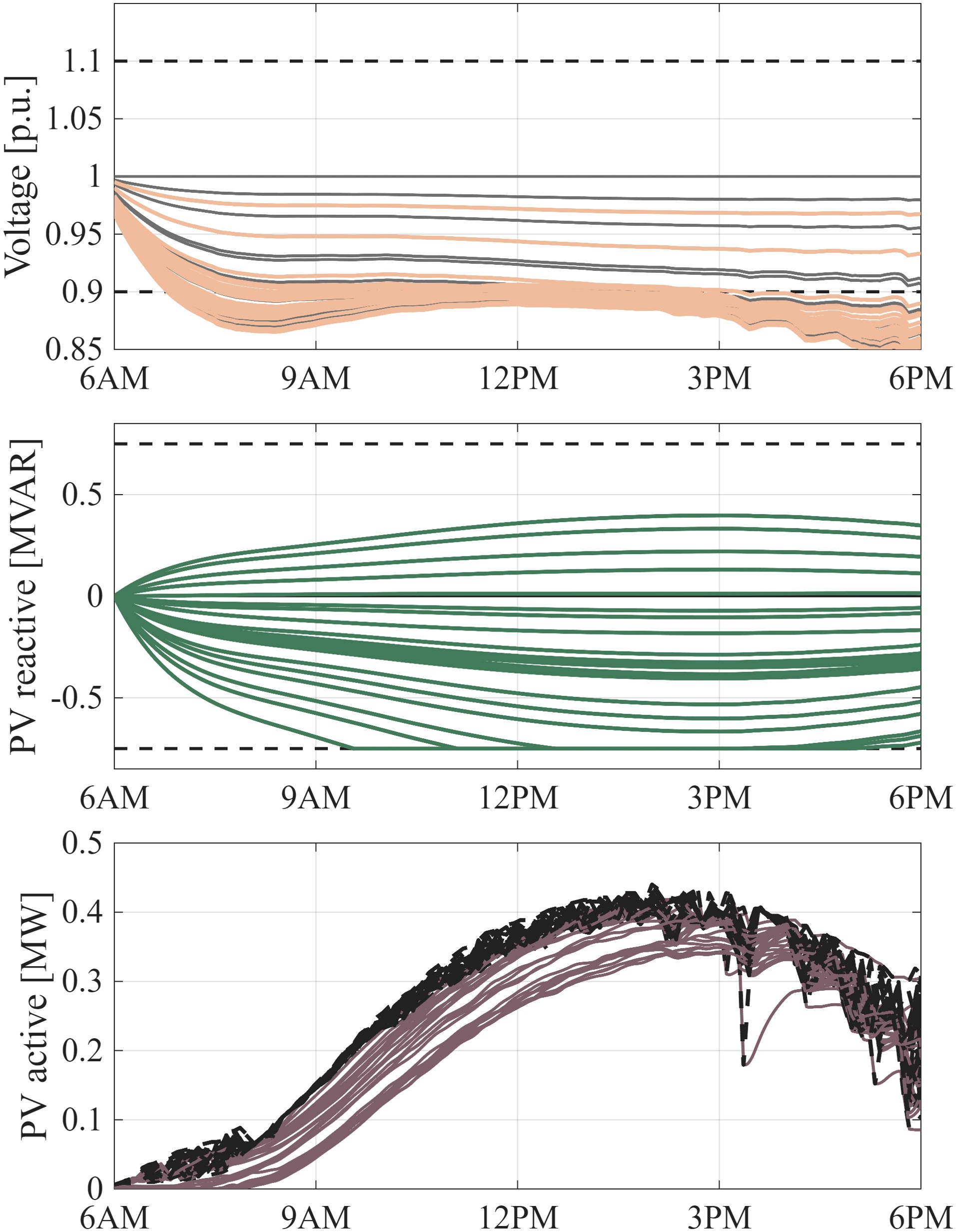}
        % \includegraphics[width=\linewidth]{graphics/topochange25-Hall-cnstr-baseFO-newTopo.png}
        % \includegraphics[width=\linewidth]{graphics/topochange-H500-uncnstr-baseFO-newTopo.png}
        % standard.pdf
        \caption{Standard feedback optimization}
        \label{fig:distribution_grid_results_standard_FO}
    \end{subfigure}%
    % \hspace{.5em}
    \hfill
    \begin{subfigure}[t]{0.32\textwidth}
        \centering
        \includegraphics[width=\linewidth]{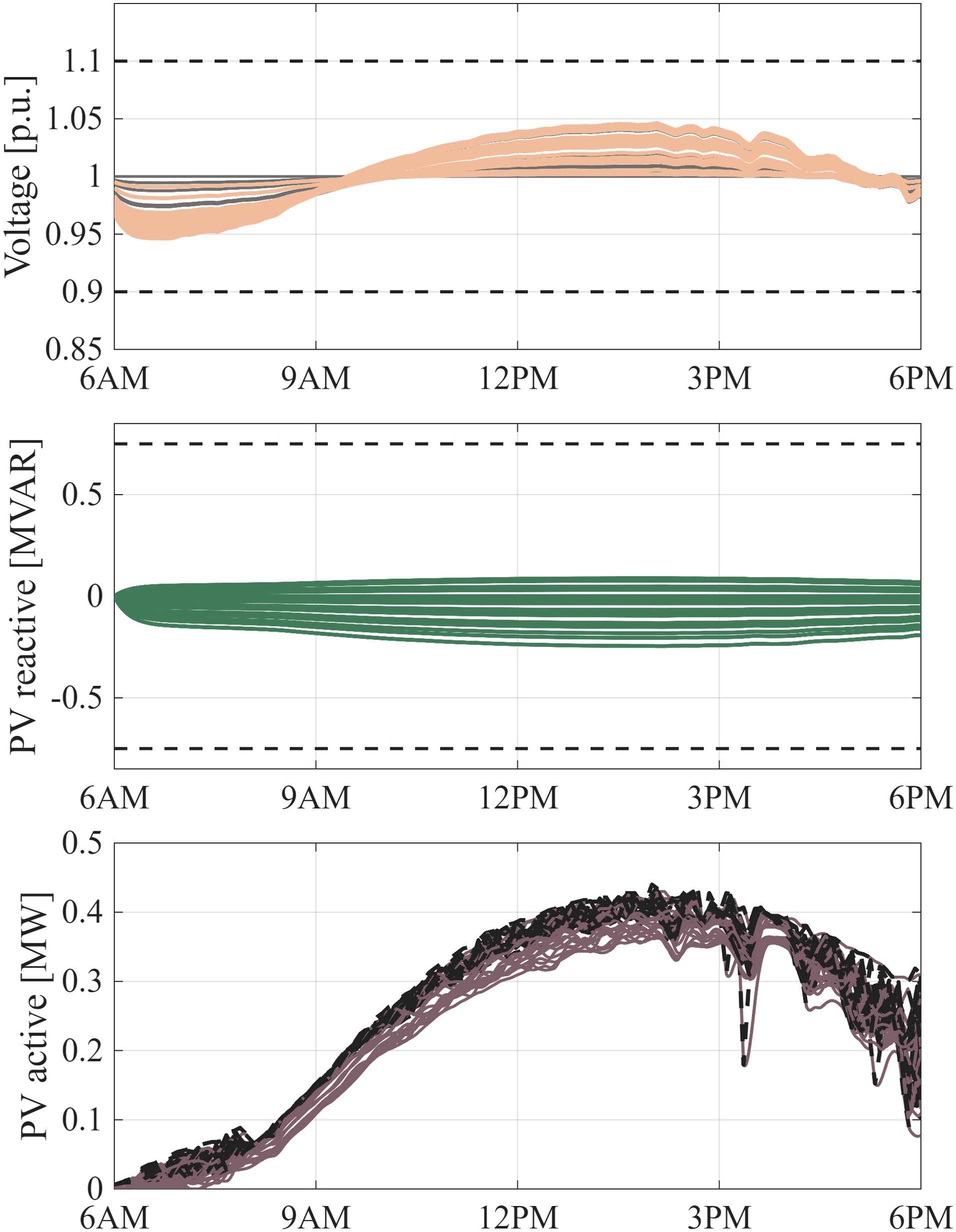}
        \caption{Robust feedback optimization \eqref{eqn:l2-step}}
        \label{fig:distribution_grid_results_robust_FO_l2}
    \end{subfigure}
    \hfill
    \begin{subfigure}[t]{0.32\textwidth}
        \centering
        \includegraphics[width=\linewidth]{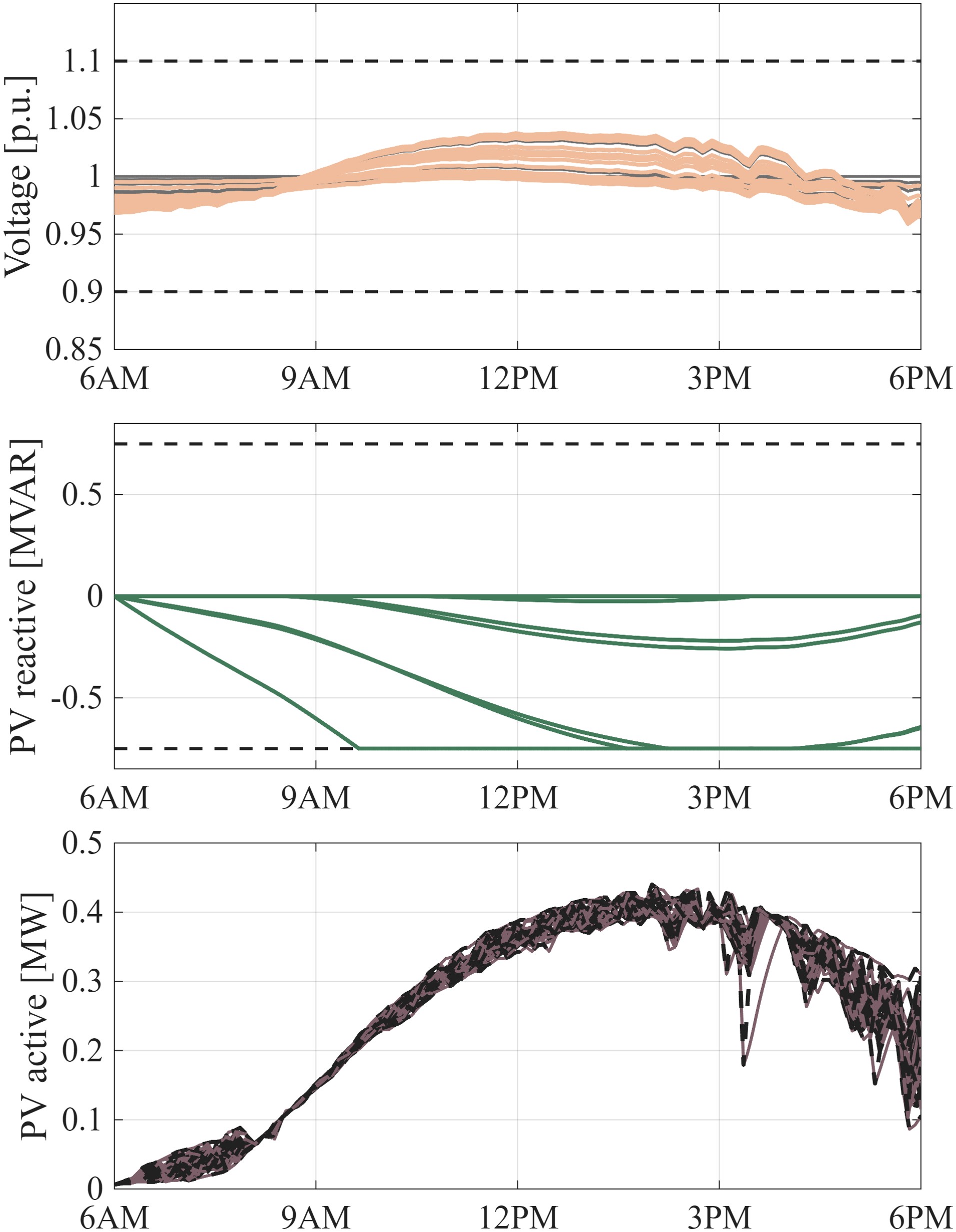}
        \caption{Robust feedback optimization \eqref{eqn:l1-step}}
        \label{fig:distribution_grid_results_robust_FO_l1}
    \end{subfigure}
    % \hspace{.5em}
    % \vspace{-.5ex}
    \caption{This figure illustrates real-time voltage control for a distribution grid after the topology change. \rev{The horizontal axis represents the time of day. The curves in the top, middle, and bottom sub-figures indicate the bus voltages, inverter reactive power injections, and inverter active power generations, respectively. The black dashed lines in the top and bottom sub-figures correspond to voltage limits and photovoltaic maximum power points, respectively.}}
    \label{fig:distribution_grid_results}
    \vspace{-.5ex}
\end{figure*}

\section{Numerical Experiments} \label{sec:simulation}
We present a case study in a distribution grid to showcase the effectiveness of our robust feedback optimization controllers. Specifically, we consider real-time voltage regulation while minimizing active power curtailment and reactive power actuation. Our goal is to show that robustification is effective beyond an academic setting for theoretical guarantees and can address practical challenges such as nonlinear steady states and state-dependent sensitivities. \rev{Our code is available at \url{https://github.com/zyhe/robustOFO}}. 

Consider a distribution grid with $n \in \mathbb{N}$ photovoltaic inverters. Let $p_{i,k}$, $q_{i,k}$, $p^{\textup{MPP}}_{i, k}$, and $v_{i,k}$ denote the active power, reactive power, maximum power point, and voltage of inverter $i$ at time $k$, respectively, where $i=1,\ldots,n$. Let $u_{i,k} \triangleq [p_{i,k} - p^{\text{MPP}}_{i, k}, q_{i,k}]$ be the variable of inverter $i$. Let $u_k = [u_{1,k}, \ldots, u_{n,k}]$ and $v_k = [v_{1,k}, \ldots, v_{n,k}]$ be concatenated variables. Further, $d_k$ represents the load at time $k$. The nonlinear map from $u_k$ and $d_k$ to $v_k$ is given by the power flow solver \cite{zimmerman2010matpower}. We aim to regulate grid voltage and minimize renewable energy curtailments and reactive power actuation. This is formalized by the following problem
\begin{equation} \label{eqn:PS-obj1}
\begin{split}
    \min_{u_k} \quad & \lVert u_k \rVert ^2_R + \lambda \lVert v_k - r_k \rVert ^2_Q \\
  \textrm{s.t.} \quad & u_{i,k} \in \mathcal{U}_{i,k} \quad \forall i = 1, \cdots, m,
  % & v_{j, \min} \leq v_{j,k} \leq v_{j, \max} & \forall \ j = 1, \cdots, n \label{eqn:PS-obj1-ycnstrs}
\end{split} 
\end{equation}
where $R \in \mathbb{R}^{2n \times 2n}$ and $Q \in \mathbb{R}^{n \times n}$ are positive definite cost matrices, $\mathcal{U}_{i,k} \triangleq \{ [p_i, \, q_i] : 0 \leq p_i \leq p^{\text{MPP}}_{i, k}, q_{\min} \leq q_i \leq q_{\max} \}$ is the constraint set, and $q_{\min}$ and $q_{\max}$ are lower and upper bounds on reactive power actuation, respectively.

We adopt the UNICORN 56-bus test case \cite{ortmann_unicorn_nodate} with 25 photovoltaic inverters. Although the input-output sensitivity is a nonlinear function of $u_k$, we learn a constant approximation $\hat{H}$ based on power flow linearization and historical data of the injected powers and voltages. This sensitivity becomes even more inexact when the grid topology changes, specifically when the point of common coupling is switched from bus 1 to bus 26. While the uncertainty set for $\hat{H}$ is hard to characterize correctly, we tune the regularizers of \eqref{eqn:l2-step} and \eqref{eqn:l1-step} by gradually decreasing their values from conservative upper bounds. We augment the standard feedback optimization controller \eqref{eq:fo_inexact_sens} and the proposed controllers \eqref{eqn:l2-step} and \eqref{eqn:l1-step} with projection to $\mathcal{U}_{i,k}$, use the same step size, and apply these controllers to the changed grid. 

As shown in the first sub-figure of \cref{fig:distribution_grid_results_standard_FO}, when implemented in a new environment with sensitivity uncertainty, standard feedback optimization causes oscillations and voltage violation. Note that the dashed lines in the sub-figures on the first row denote the maximum and minimum voltage limits, which equal $1.1$ p.u.~and $0.9$ p.u. \rev{(i.e., $1.1$ and $0.9$ times the base voltage)}, respectively. This standard controller also requires large reactive power actuation. In contrast, robust feedback optimization controllers maintain voltage stability after the point of common coupling changes. This is achieved by conservatively regulating control inputs, a consequence of regularization in the face of uncertainty. \rev{The sparsity-promoting effect of the $\ell_1$-regularizer is reflected in \cref{fig:distribution_grid_results_robust_FO_l1}, where the reactive power injections of some inverters are zero. In comparison, the $\ell_2$-regularizer induces isotropic shrinkage, see \cref{fig:distribution_grid_results_robust_FO_l2}.} Furthermore, as shown in the sub-figures in the last row, the proposed controllers lead to less active power curtailments compared to the standard approach. Overall, robust feedback optimization effectively handles model uncertainty in this example of real-time voltage regulation.

\section{Conclusion} \label{sec:conclusion}
We addressed steady-state optimization of a dynamical system subject to model uncertainty by presenting robust feedback optimization, which seeks optimal closed-loop performance given all possible sensitivities falling within bounded uncertainty sets. Tractable reformulations for this min-max problem via regularized steady-state specifications were then provided.
% and delineated how general input constraints can be handled through projection. 
We showcased the adaptation and robustness of our controllers through theoretical tracking guarantees and numerical experiments on a distribution grid. Future avenues include tuning regularizers via differentiable programming, incorporating regularization into online sensitivity learning, and pursuing robustness against model uncertainty for nonlinear stable systems.

\ifdefined\arxivVersion
\appendix
\subsection{Supporting Lemmas}
We provide two lemmas that quantify the system dynamics \eqref{eqn:plant-dynamics} and the controller iterations \eqref{eqn:l2-step} and \eqref{eqn:l1-step}.

\begin{lemma}\label{lem:dynamics_itr}
Let the conditions of \cref{thm:bounds} hold. The system dynamics \eqref{eqn:plant-dynamics} satisfy
\begin{align}\label{eq:dynamics_itr}
\lVert &x_{k+1} - x_{\textup{ss},k+1} \rVert_P \leq c_1 \lVert x_{k} - x_{\textup{ss},k} \rVert_P + \eta c_2 \|u_k - u_k^*\| \notag \\
    &+ \eta c_3 \|H \!-\! \hat{H}\|\|u_k\| + c_4 \|d_{k+1} \!-\! d_{k}\| + \eta c_5,
\end{align}
where $c_1,c_2,c_3,c_4$, and $c_5$ are constants specified in \eqref{eq:bound_constants}.
\end{lemma}
\begin{proof}
% We start by analyzing the dynamics \eqref{eqn:plant-dynamics}. 
Recall that $x_{\textup{ss},k} = (I-A)^{-1}(Bu_{k} + d_{x,k})$ is the steady state of \eqref{eqn:plant-dynamics} induced by $u_k$ and $d_{x,k}$. Let $L_x^u \triangleq \|(I-A)^{-1}B\|$ and $L_x^d \triangleq \|(I-A)^{-1}\|$ be the Lipschitz constants of $x_{\textup{ss},k}$ with respect to $u$ and $d$, respectively. 
Since $\rho(A) < 1$, for any given positive definite $\bar{Q} \in \mathbb{R}^{n \times n}$, there exists a unique positive definite $P \in \mathbb{R}^{n \times n}$ satisfying the Lyapunov equation $A^\top P A- P + \bar{Q} = 0$. Therefore,
\begin{align}\label{eq:state_evolve}
    \lVert &x_{k+1} - x_{\textup{ss},k+1} \rVert_P \stackrel{(a)}{\leq} \lVert x_{k+1} \!-\! x_{\textup{ss},k} \rVert_P + \lVert x_{\textup{ss},k+1} \!-\! x_{\textup{ss},k} \rVert_P \notag \\
        &\stackrel{(b)}{\leq} \sqrt{1 - \gamma}\|x_{k} - x_{\textup{ss},k}\|_P + \lVert x_{\textup{ss},k+1} - x_{\textup{ss},k} \rVert_P \notag \\
        &\stackrel{(c)}{\leq} \sqrt{1 - \gamma} \|x_{k} - x_{\textup{ss},k}\|_P \!+\! \lambda_{\max}(P)L_x^u \|u_{k+1} - u_{k}\| \notag \\
        &\quad + \lambda_{\max}(P)L_x^d \|d_{k+1} - d_{k}\|,
\end{align}
In \eqref{eq:state_evolve}, (a) uses the triangle inequality,
and the contraction in (b) follows from the Lyapunov equation and the property of the weighted norm, where $\gamma = \frac{\lambda_{\min}(\bar{Q})}{\lambda_{\max}(P)} \in (0,1)$, and $\lambda_{\min}(\cdot)$ and $\lambda_{\max}(\cdot)$ represent the minimum and maximum eigenvalues of a matrix, respectively. Moreover, (c) uses the triangle inequality and the Lipschitz continuity of $x_{\textup{ss},k}$. 

We proceed to analyze the term $\|u_{k+1} - u_{k}\|$ in \eqref{eq:state_evolve}. For controller~\eqref{eqn:l2-step}, let
\begin{align*}
    \hat{\nabla} \Phi_{\ell_2,k}(u_{k},y_k) &= 2R u_{k} \!+\! 2\lambda \hat{H}^{\top} Q (y_k \!-\! r_k) \!+\! 2\rho_{\text{gen}} u_{k}, \\
    \nabla \Phi_{\ell_2,k}(u_{k}) &= 2R u_{k} \!+\! 2\lambda \hat{H}^{\top} Q (\hat{H}u_{k} \!+\! d_k \!-\! r_k) \!+\! 2\rho_{\text{gen}} u_{k}
\end{align*}
be the update direction and the true gradient at $u_{k}$, respectively. Further, let $L_{\Phi'}^u = \|2R + 2\lambda\hat{H}^\top Q \hat{H}\|$ be the Lipschitz constant of $\nabla \Phi_{\ell_2,k}$ with respect to $u$. Therefore, we have
\begin{align} \label{eqn:stab-2b}
    &\lVert u_{k+1} - u_{k} \rVert = \|u_k - \eta \hat{\nabla} \Phi_{\ell_2,k}(u_{k},y_k) - u_{k} \| \notag \\
    & \stackrel{(a)}{\leq} \eta \lVert \nabla \Phi_{\ell_2,k}(u_{k}) \rVert \!+\! \eta \|\hat{\nabla} \Phi_{\ell_2,k}(u_{k},y_k) \!-\! \nabla \Phi_{\ell_2,k}(u_{k}) \| \notag \\
    & \stackrel{(b)}{\leq} \eta \lVert \nabla \Phi_{\ell_2,k}(u_{k}) - \nabla\Phi_{\ell_2,k}(u_k^*) \rVert \notag \\
    &\quad + 2\eta \lambda\|\hat{H}^\top Q\|\|y_k - (\hat{H}u_{k}+d_k)\| \notag \\
    & \stackrel{(c)}{\leq} \eta L_{\Phi'}^u \|u_{k} \!-\! u_{k}^*\| \!+\! 2\eta \lambda\|\hat{H}^\top Q\| \frac{\|C\|}{\lambda_{\min}(P)} \|x_k \!-\! x_{\textup{ss},k}\|_P, \notag \\
    &\quad \!+\! 2\eta \lambda \|\hat{H}^\top Q\| \|H - \hat{H}\|\|u_k\|,
\end{align}
where (a) uses the triangle inequality; (b) is because of the triangle inequality and $u_k^*$ being the optimal point, i.e., $\nabla\Phi_{\ell_2,k}(u_k^*) = 0$; (c) follows from the Lipschitz continuity of $\nabla \Phi_{\ell_2,k}$, the addition and subtraction of $Hu_k+d_k$ inside $\|y_k - (\hat{H}u_{k}+d_k)\|$, the expression of $y_k$, and the property of the weighted norm. For controller \eqref{eqn:l1-step}, we can perform similar analysis and obtain an upper bound akin to \eqref{eqn:stab-2b}, albeit with an additional term $2\eta \|\rho_{\textup{col}}\|$ \rev{because of the optimality condition $0 \in \nabla\Phi_{\ell_2,k}(u_k^*) + \rho_{\textup{col}}^\top \partial |u_k|$}. We incorporate the above results into \eqref{eq:state_evolve} and obtain \eqref{eq:dynamics_itr}, where the constants are given by
\begin{equation}\label{eq:bound_constants}
\begin{split}
    c_1 &= \sqrt{1-\gamma} + 2\eta \lambda L_x^u \|\hat{H}^\top Q\| \|C\| \frac{\lambda_{\max}(P)}{\lambda_{\min}(P)},\\
    c_2 &= \lambda_{\max}(P)L_x^u L_{\Phi'}^u, \\
    c_3 &= 2 \lambda L_x^u \lambda_{\max}(P) \|\hat{H}^\top Q\|, \\
    c_4 &= \lambda_{\max}(P)L_x^d, \\
    c_5 &= \begin{cases} 0, & \textup{for } \eqref{eqn:l2-step}, \\ 2 \lambda_{\max}(P)L_x^u\|\rho_{\textup{col}}\|, & \textup{for } \eqref{eqn:l1-step}. \end{cases}
\end{split}
\end{equation}
Therefore, \cref{lem:dynamics_itr} is proved.
 % (see \eqref{eqn:plant-dynamics})
\end{proof}

In the following lemma, we characterize the property of the controller iterations \eqref{eqn:l2-step} and \eqref{eqn:l1-step}. Both the objective \eqref{eqn:obj-l2} and the quadratic part of \eqref{eqn:obj-l1} are strongly convex and smooth in $u$. Let $\mu_{\Phi}$ and $L_{\Phi}$ be the corresponding constants of strong convexity and smoothness, respectively. Recall that $P$ is the matrix appearing in the weighted norm in \cref{lem:dynamics_itr}.

\begin{lemma}\label{lem:controller_itr}
    Let the conditions of \cref{thm:bounds} hold. The controller iterations \eqref{eqn:l2-step} and \eqref{eqn:l1-step} satisfy
    \begin{equation}\label{eq:controller_itr}
    \begin{split}
        \|u_{k+1} \!-\! u_{k+1}^*\| \leq& \alpha \|u_{k} \!-\! u_{k}^*\| \!+\! \eta \frac{L_T^y \|C\|}{\lambda_{\min}(P)} \rev{\lVert x_k \!-\! x_{\textup{ss},k}\rVert_P} \\
        &\!+\! \eta L_T^y \|H\!-\!\hat{H}\| \|u_{k}\| \!+\! \lVert u_{k+1}^* \!-\! u_k^* \rVert,
    \end{split}
    \end{equation}
    where $\alpha = \sqrt{1-\eta(2\mu_{\Phi} - \eta L_{\Phi}^2)}$, and $L_T^y = 2\lambda \|\hat{H}^\top Q\|$.
\end{lemma}
\begin{proof}
% We proceed to analyze the controller iteration. 
Let the right-hand side of \eqref{eqn:l2-step} or \eqref{eqn:l1-step} be denoted by $T(u_k,y_k)$. The mapping $T(u,y)$ is $\eta L_T^y$-Lipschitz in $y$, where $L_T^y = 2\lambda \|\hat{H}^\top Q\|$. Therefore, 
% the controller iterations satisfy
\begin{equation*}
\begin{split}
    \lVert u_{k+1} &- u_{k+1}^* \rVert \stackrel{(a)}{\leq} \lVert T\bigl(u_k, y_k\bigr) - u_k^* \rVert + \|u_{k+1}^* - u_k^*\| \\
    & \stackrel{(b)}{\le} \lVert T\bigl(u_{k}, y_k\bigr) - T\bigl(u_{k}, Hu_{k} + d_k\bigr) \rVert \\
    &\quad + \lVert T\bigl(u_{k}, Hu_{k} \!+\! d_k\bigr) \!-\! T\bigl(u_{k}, \hat{H}u_{k} \!+\! d_k\bigr) \rVert \\
    &\quad + \lVert T\bigl(u_{k}, \hat{H}u_{k} + d_k\bigr) - u_k^* \rVert + \|u_{k+1}^* - u_k^*\| \\
    & \stackrel{(c)}{\le} \eta L_T^y \lVert y_k - (Hu_{k} + d_k)\rVert + \eta L_T^y \|H-\hat{H}\| \|u_{k}\| \\
    &\quad + \alpha \|u_k - u_k^*\| + \lVert u_{k+1}^* - u_k^* \rVert \\
    & \stackrel{(d)}{\le} \eta \frac{L_T^y \|C\|}{\lambda_{\min}(P)} \rev{\lVert x_k - x_{\textup{ss},k}\rVert_P} + \eta L_T^y \|H-\hat{H}\| \|u_{k}\| \\
    &\quad+ \alpha \lVert u_{k} - u_{k}^* \rVert + \lVert u_{k+1}^* - u_k^* \rVert, \\
\end{split}
\end{equation*}
where (a) and (b) use the triangle inequality; (c) follows from the Lipschitz continuity of $T$ and the contraction of $T$ (see \cite[Proposition 25.9]{bauschke_convex_2011}, where $\alpha$ is given in the lemma when $\eta \in (0, 2\mu_\Phi / L_\Phi^2)$, and we also use the non-expansiveness property of the proximal operator for \eqref{eqn:l1-step}); and (d) applies \eqref{eqn:plant-dynamics} and the property of the weighted norm.
\end{proof}

\subsection{Proof of Theorem~\ref{thm:bounds}}\label{app:proof_thm}
The main idea is to analyze the coupled evolution of state dynamics and controller iterations, whose properties are established in \cref{lem:dynamics_itr,lem:controller_itr}, respectively.

\begin{proof}
The coupling between the state dynamics and controller iterations can be compactly written as
\begin{align}\label{eq:contraction_metric}
    &\underbrace{\begin{bmatrix}
        \|u_{k+1} \!-\! u_{k+1}^*\| \\ \lVert x_{k+1} \!-\! x_{\textup{ss},k+1} \rVert_P
    \end{bmatrix}}_{\triangleq w_{k+1}}
    \leq
    \underbrace{
    \begin{bmatrix}
        \alpha & \eta \frac{L_T^y \|C\|}{\lambda_{\min}(P)} \\
        \eta c_2 & c_1
    \end{bmatrix}}_{\triangleq M}
    \underbrace{
    \begin{bmatrix}
        \|u_k \!-\! u_k^*\| \\ \lVert x_k \!-\! x_{\textup{ss},k} \rVert_P
    \end{bmatrix}}_{\triangleq w_k} \notag \\
    &+
    \underbrace{
    \begin{bmatrix}
        \eta L_T^y \|H\!-\!\hat{H}\| \|u_{k}\| \!+\! \lVert u_{k+1}^* \!-\! u_k^* \rVert \\
        \eta c_3 \|H \!-\! \hat{H}\|\|u_k\| + c_4 \|d_{k+1} \!-\! d_{k}\| + \eta c_5
    \end{bmatrix}}_{\triangleq q_k},
\end{align}
where the constants $c_1$ to $c_4$ are given by \eqref{eq:bound_constants}.
Note that $M$ in \eqref{eq:contraction_metric} is a $2$-by-$2$ positive matrix, and therefore its Perron eigenvalue equals $\rho(M)$. Hence, the requirement that $\rho(M) < 1$ is equivalent to $m_{11} + m_{22} - m_{11}m_{22} + m_{21}m_{12} < 1$, where $m_{ij}$ denotes the $ij$-th element of $M$. This inequality translates to
\begin{equation}\label{eq:sz_cond}
    g(\eta) \triangleq \frac{L_T^y \|C\| c_2}{\lambda_{\min}(P)} \eta^2 + \alpha + c_1 - \alpha c_1 < 1,
\end{equation}
where $\alpha$ and $c_1$ are given in \cref{lem:controller_itr} and \eqref{eq:bound_constants}, respectively. When $\eta=0$, we have $\alpha=1, c_1=\sqrt{1-\gamma}$. The function $g(\eta)$ satisfies $g(0) = 1, g'(0) = -(1-\sqrt{1-\gamma})\mu_{\Phi} < 0$. Hence, there exists $\eta^* \in (0, 2\mu_\Phi / L_\Phi^2)$ such that for any $\eta \in (0, \eta^*)$, $g(\eta) < 1$, implying $\rho(M) < 1$.
We telescope \eqref{eq:contraction_metric} and obtain
\begin{equation} \label{eqn:optstab-3}
    w_k \leq M^k w_0 + \sum^{k-1}_{i=0} M^{k-i} q_{i+1}.
\end{equation}
When $\eta \in (0, \eta^*)$, there exists $r > 0$ and $c_M \in [0, 1)$ such that $\lVert M^k \rVert \leq r (c_M)^k$, see also \cite[Chapter 5]{lasalle_stability_1986}. Hence, we obtain from \eqref{eqn:optstab-3} the following inequality
% This result is exploited by taking the norm of (\ref{eqn:optstab-3}) to give (a) below:
\begin{equation*} %\label{eqn:optstab-4}
\begin{split}
    \lVert w_k \rVert & \leq r (c_M)^k \lVert w_0 \rVert + \sum^{k-1}_{i=0} r (c_M)^{k-i} \lVert q_{i+1} \rVert \\
    & \stackrel{(a)}{\le} r (c_M)^k \lVert w_0 \rVert + r c_M \lVert \bar{q} \rVert \sum^{k-1}_{i=0} (c_M)^i \\
    & \stackrel{(b)}{\le} r (c_M)^k \lVert w_0 \rVert + r \frac{c_M}{1 - c_M} \lVert \bar{q} \rVert, \\
    % & \stackrel{(c)}{\le} r (c_M)^k \sqrt{\lambda_{\max}(P)} \left \lVert
    % \begin{bmatrix} \mathbf{d} u_0 \\ \mathbf{d} x_0 \end{bmatrix}
    % \right \rVert
    % + r \frac{c_M}{1 - c_M} \lVert \bar{q} \rVert,
\end{split}
\end{equation*}
while (a) is due to $\bar{q} \triangleq \sup_{i \in [k]} q_i$, and (b) uses the upper bound on the partial sum of a geometric series. Therefore, \eqref{eqn:optstab-1} is proved.
% (b) is obtained by defining $\bar{q} := $ and re-expressing the summation in the second term; (c) by observing that a geometric series upper bounds the first $k$ terms of the re-expressed summation; and (d) by the relation
% \begin{equation} \label{eqn:Pnormbounds}
%     \sqrt{\lambda_{\min}(P)} \lVert \cdot \rVert \leq \lVert \cdot \rVert _P \leq \sqrt{\lambda_{\max}(P)} \lVert \cdot \rVert
% \end{equation}
% applied to $\begin{bmatrix} \mathbf{d} u_0 \\ \mathbf{d} x_0 \end{bmatrix}$, where we recall $\bar{\omega}_k := \begin{bmatrix} \lVert \mathbf{d} u_k \rVert _P \\ \lVert \mathbf{d} x_k \rVert _P \end{bmatrix}$.
% Finally, the lower bound of (\ref{eqn:Pnormbounds}) is used on (\ref{eqn:optstab-4}) to obtain (\ref{eqn:optstab-1}), with $r_1 := r \frac{\sqrt{\lambda_{\max}(P)}}{\sqrt{\lambda_{\min}(P)}}$ and $r_2 := \frac{r}{\sqrt{\lambda_{min}(P)}}$.
\end{proof}
\fi

% \ifCLASSOPTIONcaptionsoff
%   \newpage
% \fi

\section*{Acknowledgement}
We thank Prof.~Linbin Huang for inspirational discussions.

%%% bibliography %%%
\balance

\bibliographystyle{IEEEtran}
\bibliography{robustFO}

% that's all, folks!
\end{document}